# Binary Weight Allocation for Multi-Objective Path Optimization: Efficient Earliest and Latest Path Discovery in Network Systems

Wei-Chang Yeh
Integration and Collaboration Laboratory
Department of Industrial Engineering and Engineering Management
National Tsing Hua University
Department of Industrial and Systems Engineering, Chung Yuan Christian University
yeh@ieee.org

**Abstract:** This paper proposes earliest and latest path algorithms based on binary weight allocation, assigning weights of $2^{(i-1)}$ and $2^{(m-i)}$ to the *i*-th arc in a network. While traditional shortest path algorithms optimize only distance, our approach leverages Binary-Addition-Tree ordering to efficiently identify lexicographically smallest and largest paths that establish connectivity. These paths partition the solution space into three regions: guaranteed disconnection, transitional connectivity, and guaranteed no simple paths. Our weight allocation enables implicit encoding of multiple objectives directly in binary representations, maintaining the $O((|V|+|E|)\log|V|)$ complexity of Dijkstra's algorithm while allowing simultaneous optimization of competing factors like reliability and cost. Experimental validation demonstrates significant computational time reduction compared to traditional multi-objective methods. Applications span telecommunications, transportation networks, and supply chain management, providing efficient tools for network planning and reliability analysis under multiple constraints.

## 1. INTRODUCTION

With the rapid development of advanced technologies, various networks have emerged and become more diverse and powerful. These networks include the Internet of Things [1], 4G/5G telecommunications [2, 3], social networks [4, 5], deep learning [6], cloud/fog/edge computing [7, 8], and smart wireless sensor networks [9]. Modern networks and traditional networks (such as water, gas, electricity, and telephone networks) have become an indispensable part of the daily lives of almost all humans and various industries (such as manufacturing, commerce, and supply chains) globally [10, 11].

In the analysis and optimization of complex network systems, the shortest path algorithm has always been a core problem in graph theory research and a basic tool for network planning, resource allocation, and reliability assessment. Consider an undirected graph $G(V, E)$, where $V = \{v_1, v_2, ..., v_n\}$ represents the set of nodes, and $E = \{a_1, a_2, ..., a_m\}$ represents the set of arcs. Each arc $a_i$ connects two nodes and has a weight $W(a_i)$. The traditional shortest path problem is usually defined as: given a source node $1 \in V$ and a sink node $n \in V$, find a path P from 1 to $n$ such that the sum of weights of all arcs on the path $\sum_{\forall a \in P} W(a)$ is minimized.

Traditional shortest path algorithms (such as Dijkstra's algorithm [12], Bellman-Ford algorithm [13, 14], etc.) mainly focus on optimizing the single dimension of path length. Although they perform excellently in many scenarios, they often fall short when facing increasingly complex multi-objective decision-making needs in modern network systems. As network systems grow more complex, so do the requirements for path optimization algorithms. While traditional algorithms excel in minimizing single objectives like distance, modern applications require simultaneous consideration of multiple factors such as path distance, arc priority, and system reliability, where traditional approaches struggle to provide satisfactory solutions.

We leverage BATs [4, 15, 16, 17, 18, 19, 20, 21, 22, 23, 24, 25], a structured approach to systematically enumerate binary vectors in increasing order, where each node in the tree represents a unique binary vector corresponding to a specific subset of network arcs. To address these multi-objective challenges, we introduce earliest and latest shortest path algorithms leveraging a novel binary weight allocation strategy in the BAT ordering [24]. This ordering is precisely defined such that vector $X$ precedes vector $X^*$ if $X \ll X^*$ without considering the weight, meaning there exists an index $i$ where $X(a_i) < X^*(a_i) = 1$ and $X(a_j) = X^*(a_j)$ for all $j > i$ [24]. This lexicographical ordering creates a systematic exploration of the solution space that prioritizes paths with specific arc inclusion patterns, enabling more efficient connectivity analysis and constraint satisfaction.

For example, consider a network with 5 arcs $\{a_1, a_2, a_3, a_4, a_5\}$. Let vector $X = (1, 0, 0, 1, 0)$ and $X^* = (1, 0, 1, 1, 0)$, where each position indicates whether the corresponding arc is included (1) or excluded (0) in the solution. According to the BAT ordering: we compare positions from right to left;

for positions $j = 5, 4$: $X(a_j) = X^*(a_j)$, so we continue; at position $i = 3$: $X(a_3) = 0 < X^*(a_3) = 1$; therefore, $X \ll X^*$, meaning $X$ found precedes $X^*$ in the BAT ordering, i.e., $X$ is found before $X^*$ using the BAT.

Any vector $X$ found before the earliest path is disconnected in $G(X) = G(V, \{a \in E \mid X(a) = 1\})$, while any vector $X$ found after the latest path is always connected in $G(X)$. The identification of earliest and latest paths within the BAT serves several crucial purposes in network optimization and routing applications.

This paper makes several key contributions to the field of network optimization: (1) We introduce a novel binary weight allocation strategy that enables implicit encoding of path priorities; (2) We formally define and prove properties of earliest and latest paths within the BAT framework; (3) We develop efficient algorithms that leverage these properties to significantly reduce computational overhead; and (4) We demonstrate practical applications across multiple domains including supply chain management and power system reliability assessment.

For the earliest path, determining where connectivity first emerges provides a threshold that eliminates the need to examine disconnected vectors, significantly reducing computational overhead. This boundary marks the transition from guaranteed disconnection to potential connection, making it valuable for:

1. **Feasibility testing**: Quickly identifying the minimum resources needed to establish connectivity in resource-constrained networks [16].
2. **Network resilience analysis**: Finding the critical threshold where network connectivity becomes possible, essential for fragile or emergency networks [19].
3. **Incremental network design**: Determining the minimal set of links required to achieve basic connectivity [25].

For the latest path, identifying where guaranteed no simple path begins provides another important computational boundary. This threshold marks where all subsequent vectors will have connectivity, offering advantages for:

1. **Robust routing guarantees**: Ensuring reliable transmission in mission-critical systems by operating beyond this threshold

2. **Service level agreement (SLA) enforcement**: Providing mathematical guarantees for network availability
3. **Simplified route planning**: Eliminating connectivity checking for all vectors beyond this threshold

For instance, in telecommunications network design, our approach allows engineers to prioritize high-reliability links while maintaining acceptable latency. By assigning lower indices to higher-reliability links, the earliest path algorithm naturally discovers routes that maximize reliability for a given distance constraint without requiring complex multi-objective optimization frameworks.

Together, these two path boundaries effectively partition the solution space into three regions: guaranteed disconnection (before earliest), transitional connectivity (between earliest and latest), and guaranteed no simple paths (after latest). This partitioning allows algorithms to optimize their approach based on which region they're operating in, potentially skipping unnecessary connectivity tests or employing different computational strategies.

Our empirical evaluations demonstrate that the proposed approach achieves $O(n)$ computational time compared to traditional multi-objective optimization methods. This efficiency gain becomes particularly significant as network complexity increases, with time complexity growing nearly linearly with network size in contrast to the polynomial growth observed in conventional approaches.

These interdisciplinary application cases jointly validate the dual value of the earliest path algorithm based on binary weights at both theoretical and practical levels. In terms of computational efficiency, the algorithm exhibits near-linear time complexity growth in large-scale networks, significantly outperforming traditional multi-objective optimization methods. In terms of application flexibility, the algorithm can adapt to diverse needs in different domains by simply adjusting the arc index allocation strategy. In terms of optimization effect, the algorithm can effectively balance multiple constraint factors while ensuring basic path performance, providing high-quality solutions for decision-makers.

The paper is structured as follows: Section 2 establishes crucial notations, acronyms, and underlying assumptions that form the basis of our work. Section 3 examines traditional shortest path algorithms and identifies their constraints. Section 4 develops the theoretical framework for BATs and proves essential properties of first connected vectors. Section 5 presents our innovative earliest path algorithm and demonstrates its enhanced efficiency compared to current methods. Section 6 details the implementation with optimized pseudocode for fast BAT construction, provides practical examples, and presents experimental results confirming our performance improvements. Section 7 concludes by highlighting our contributions, discussing practical applications, and suggesting future research directions.

## 2. FOUNDATIONAL ELEMENTS: NOTATION, TERMINOLOGY, AND PREMISES

This section introduces the core acronyms, mathematical notation, foundational assumptions, and terminology necessary for understanding and implementing our proposed methodology.

### 2.1 Acronyms

BAT: Binary-Addition-Tree Algorithm (efficient combinatorial search for network states)

### 2.2 Notations

$|\bullet|$: Cardinality of a set

$||\bullet||$: Dimension of a vector or subvector

$a_i$: The arc $i$ in the network

$e_{i,j}$: Arc connecting node $i$ to node $j$

$V$: Set of nodes, where $V = \{1, 2, \ldots, n\}$

$E$: Set of arcs, where $E = \{a_1, a_2, \ldots, a_m\}$

$P$: Path connecting source node 1 to sink node $n$.

$X$: Binary state vector $X = (x_1, x_2, ..., x_m)$ with $x_i = X(a_i) \in \{0, 1\}$

$X(a_i)$: Binary state of arc $a_i$ in vector $X$ (Example: in $X = (1, 1, 1, 0, 0, 0, 1, 0)$, $X(a_i) = 1$ for $i = 1, 2, 3, 7$ and $X(a_i) = 0$ for $i = 4, 5, 6, 8$)

$W(a_i)$: Weight assigned to arc $a_i$

$G(V, E)$: Undirected graph with $V$ and $E$. Example: **Figure 1** illustrates a graph where $V = \{1, 2, 3, 4, 5\}$, $E = \{a_1, a_2, \ldots, a_8\}$, with node 1 as source and node 5 as sink.

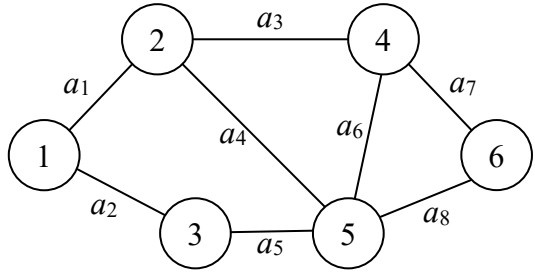

**Figure 1.** Example Network.

dist[v]: Array storing the shortest distance from source node 1 to node $v$.

prev[v]: Array storing the predecessor node of $v$ in the shortest path.

Q: Priority queue employed in Dijkstra's algorithm, ordered by dist values.

alt Temporary variable used for distance updates during Dijkstra's relaxation step.

$G(V, E, W)$: Binary-state network with structure $G(V, E)$ and weight $W$.

$G(X)$: Subgraph corresponding to state vector $X$, defined as $G(X) = G(V, \{a \in E \mid X(a) = 1\})$.

Example: For $X = (1, 1, 1, 0, 0)$, $G(X)$ displays functioning arcs (in solid lines) as shown in **Figure 2**.

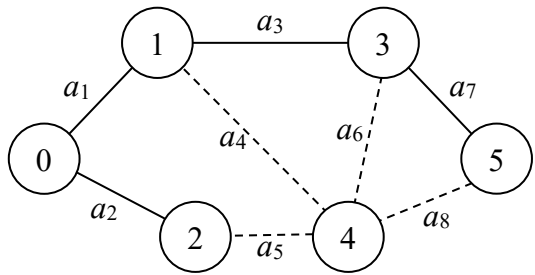

**Figure 2.** $X = (1, 1, 1, 0, 0, 0, 1, 0)$ and its corresponding $G(X)$, where the original graph $G(V, E)$ is depicted in Figure 1.

<<: Binary relation between vectors: $X << X^*$ if there exists an index $i$ where $X(a_i) = 0 < X^*(a_i) = 1$ and $X(a_j) = X^*(a_j)$ for all $j > i$.

## 2.3 Nomenclature

Binary-state network: A network where each arc exists in one of two states: functional (1) or failed (0).

Connected vector: A state vector $X$ where the network graph $G(X)$ contains at least one operational path between the source and target nodes.

Disconnected vector: A state vector $X$ where no operational path exists between the source and target nodes in $G(X)$.

BAT ordering: The BAT framework generates binary state vectors $X = (x_1, x_2, ..., x_m)$, where $x_i \in \{0, 1\}$ for $i = 1, 2, …, m$, through lexicographical ordering. For vectors $X$ and $X^*$, $X << X^*$ if $X(a_i) = 0 < X^*(a_i) = 1$ and $X(a_j) = X^*(a_j)$ for all $j > i$. This ordering prioritizes arc inclusion patterns, enabling efficient connectivity analysis.

Earliest Path: The first vector $X$ in BAT order where $G(X) = G(V, \{a \in E \mid X(a)=1\})$ becomes connected.

Latest Path: The last vector $X$ in BAT order where transitions to guaranteed no simple paths.

$X_{FC}$ (First Connected Vector): It is called the earliest path here.

## 2.4 Assumptions

This study adopts the following foundational assumptions to ensure analytical tractability within the framework of a binary-state network:

1. Network Structure: No parallel arcs and self-loops are present in the network.
2. Node Properties: All nodes are interconnected within the network
3. Arc Characteristics: Each arc operates in a binary state (working/failed).

## 3. RELATED RESEARCH ON SHORTEST PATH ALGORITHMS

Dijkstra's algorithm [12] is one of the most famous single-source shortest path algorithms, widely applied in numerous fields such as network routing, traffic planning, and resource allocation due to its efficiency and intuitiveness. The algorithm is based on a greedy strategy, starting from the source node and continuously selecting the node with the smallest current distance from a priority queue for expansion, gradually building a shortest path tree from the source node to all reachable nodes in the graph.

The core idea of Dijkstra's algorithm is the principle of "optimal substructure," meaning that any subpath of a shortest path is also the shortest [12]. When implemented with a binary heap, its time complexity is $O((|V|+|E|)\log|V|)$ [12], and with a Fibonacci heap, it can be further optimized to $O(|E|+|V|\log|V|)$ [26], making it an ideal choice for handling large-scale sparse graphs. Its space complexity is $O(|V|)$ to store the distance and predecessor arrays. However, a significant limitation of Dijkstra's algorithm is that it is only applicable to graphs with non-negative arc weights, as the algorithm assumes that once a node is processed, its shortest distance will not change again.

The Bellman-Ford algorithm [13] overcomes this limitation of Dijkstra's algorithm and can handle graphs containing negative weight arcs, and even detect whether there exists a negative weight cycle in the graph. The algorithm adopts a dynamic programming approach, gradually optimizing the shortest path estimates from the source node to all other nodes through at most $|V|-1$ relaxation operations (where $|V|$ is the number of nodes) on all arcs. Although the Bellman-Ford algorithm has wider applicability, its $O(|V||E|)$ time complexity and $O(|V|)$ space complexity are significantly higher than Dijkstra's algorithm, limiting its application efficiency on large-scale graphs [17]. In practical implementations, performance can be improved in some cases through early termination techniques (ending early when a round of relaxation operations fails to improve any shortest path estimate) and queue optimization, but for dense graphs, its efficiency remains significantly lower than Dijkstra's algorithm.

The Floyd-Warshall algorithm [14] takes a completely different approach, calculating the shortest paths between any two points in the graph with a time complexity of $O(|V|^3)$ and space complexity of $O(|V|^2)$, suitable for global path analysis and calculating the shortest distances between all pairs of nodes [18]. The algorithm is based on dynamic programming principles, gradually considering all possible paths using each node as an intermediate point through three nested loops, ultimately constructing a complete shortest path matrix [19]. This approach is particularly valuable in telecommunications and network reliability assessment, where all-pairs shortest paths are frequently required.

The main advantage of the Floyd-Warshall algorithm lies in its concise implementation and ability to process all node pairs, making it particularly suitable for relatively small dense graphs and application scenarios requiring frequent queries of shortest paths between different node pairs [14]. Additionally, similar to the Bellman-Ford algorithm [13], the Floyd-Warshall algorithm can also handle negative weight arcs and can be used to detect the existence of negative weight cycles [14]. However, for large-scale sparse graphs, especially when only single-source shortest path problems need to be solved, its cubic time complexity makes its efficiency far lower than Dijkstra's algorithm and the Bellman-Ford algorithm [13].

When extending these algorithms to multi-objective optimization scenarios, researchers have previously proposed various approaches such as the weighted sum method [21], where multiple objectives are combined into a single objective using predefined weights, and the constraint method [28], which optimizes one objective while constraining others. However, these extensions often struggle with defining appropriate weights and handling trade-offs between competing objectives. Pareto-optimal approaches [8, 11, 27] identify sets of non-dominated solutions but typically incur exponential computational costs as the number of objectives increases, making them impractical for large-scale networks.

These three classic shortest path algorithms each have their advantages and applicable scenarios, forming the basic tools for graph theory research and network analysis. The earliest/latest path algorithm based on binary weights proposed in this paper can be viewed as an innovative extension of Dijkstra's algorithm, enhancing its application capability in multi-objective decision-making through its novel binary weight allocation mechanism. Unlike previous extensions, our approach implicitly encodes multiple objectives within the binary weight structure, preserving the computational efficiency of Dijkstra's algorithm while providing systematic exploration of the solution space. This is particularly valuable in applications such as supply chain route optimization and power grid reliability analysis, where multiple competing objectives must be balanced efficiently.

## 4. BAT AND FIRST CONNECTED VECTOR

Our proposed earliest path algorithm is a new shortest path algorithm, with the original concept of the earliest path coming from the first connected vector of the (fast) BAT. Therefore, in this section, we introduce the shortest path, BAT, and the first connected vector.

### 4.1 BAT

The Binary Addition Tree (BAT) proposed by Yeh [15] adopts a new heuristic search method that is more efficient than depth-first search (DFS). The BAT algorithm is also more effective than breadth-first search (BFS) [16] and universal generating function method (UGFM) [17], as the latter produces incorrect results due to computer memory overflow. Additionally, BAT is easy to understand, easy to code, and customizable [15-20].

In a binary vector, each coordinate value is 0 or 1. Assume the search process uses a program that iteratively moves from the last coordinate to the first coordinate to update the binary vector. At the beginning of the update process, this binary vector is called a "zero vector." The core idea is to efficiently generate all possible $k$-tuple binary state vectors through the following two simple rules:

1. Rule 1: Find the first coordinate position with a value of 0 (denoted as $x_i$), and set all coordinates $j < i$ to 0;
2. Rule 2: If $i$ is the last coordinate, terminate the algorithm; at this point, all vectors have been generated.

Based on the above rules, the pseudocode for the BAT algorithm is as follows [15]:

**Algorithm: BAT**

**Input:** $k$

**Output:** All $k$-tuple binary state vectors.

**STEP 0:** Initialize vector $X = (x_1, x_2, \ldots, x_k)$ to zero vector and set $i = 1$.

**STEP 1:** If $x_i = 0$, set $x_i = 1$, reset $i = 1$, and jump to **STEP 1**.

**STEP 2:** If $i$ is the last coordinate (i.e., $i = k$), terminate the algorithm.

**STEP 3:** Set $x_i = 0$, increment $i = i + 1$, and jump to **STEP 1**.

The BAT requires only four steps and one dynamically updated k-tuple vector X, featuring concise code, efficient operation, low memory usage, and flexible adaptation. Its time complexity is $O(2^k)$, which is optimal for generating all $2^k$ possible binary vectors, while maintaining only $O(k)$ space complexity, making it significantly more memory-efficient than recursive approaches like DFS or BFS which require $O(k \cdot 2^k)$ space in the worst case.

A key property of BAT that makes it particularly suitable for our earliest/latest path algorithm is the specific ordering in which it generates vectors. The BAT produces vectors in an order where coordinate changes occur systematically from right to left, creating a lexicographical ordering that preserves structural relationships between consecutive vectors. This ordering property is crucial for our approach, as it allows us to efficiently identify connectivity transitions in the network without evaluating all possible states.

The BAT has been widely applied in fields such as network reliability calculation [21], system resilience assessment [19], wildfire spread probability analysis [22], and computer virus spread modeling [23]. Taking a 5-tuple binary state vector as an example, Table 1 shows all vectors generated by the BAT:

**Table 1.** 5-tuple binary state vectors generated based on BAT

| $i$ | $X_i$ | $i$ | $X_i$ |
|---|---|---|---|
| 1 | (0, 0, 0, 0, 0) | 17 | (0, 0, 0, 0, 1) |
| 2 | (1, 0, 0, 0, 0) | 18 | (1, 0, 0, 0, 1) |
| 3 | (0, 1, 0, 0, 0) | 19 | (0, 1, 0, 0, 1) |
| 4 | (1, 1, 0, 0, 0) | 20 | (1, 1, 0, 0, 1) |
| 5 | (0, 0, 1, 0, 0) | 21 | (0, 0, 1, 0, 1) |
| 6 | (1, 0, 1, 0, 0) | 22 | (1, 0, 1, 0, 1) |
| 7 | (0, 1, 1, 0, 0) | 23 | (0, 1, 1, 0, 1) |
| 8 | (1, 1, 1, 0, 0) | 24 | (1, 1, 1, 0, 1) |
| 9 | (0, 0, 0, 1, 0) | 25 | (0, 0, 0, 1, 1) |
| 10 | (1, 0, 0, 1, 0) | 26 | (1, 0, 0, 1, 1) |
| 11 | (0, 1, 0, 1, 0) | 27 | (0, 1, 0, 1, 1) |
| 12 | (1, 1, 0, 1, 0) | 28 | (1, 1, 0, 1, 1) |
| 13 | (0, 0, 1, 1, 0) | 29 | (0, 0, 1, 1, 1) |
| 14 | (1, 0, 1, 1, 0) | 30 | (1, 0, 1, 1, 1) |
| 15 | (0, 1, 1, 1, 0) | 31 | (0, 1, 1, 1, 1) |
| 16 | (1, 1, 1, 1, 0) | 32 | (1, 1, 1, 1, 1) |

While BAT efficiently generates all possible binary vectors, examining every vector becomes computationally infeasible for large networks. Our contribution extends BAT by identifying critical transition points – the earliest and latest paths – that partition the vector space into regions with guaranteed properties, eliminating the need to evaluate all vectors. This approach maintains BAT's systematic exploration advantage while overcoming its exponential scaling limitation for practical applications.

The BAT algorithm systematically generates all possible binary state vectors, providing an exhaustive but efficient solution for complex network analysis. The vectors it generates can be directly used to evaluate network connectivity, calculate system reliability, or simulate risk propagation processes, making it an important tool for theoretical research and engineering practice.

### 4.2 First Connected Vector

The first connected vector $X_{FC}$ is an important concept in the BAT [24], representing the first state vector in all vectors arranged in lexicographical order that makes the network connected from source node 1 to sink node $n$. Formally, $X_{FC}$ is the lexicographically smallest vector $X$ such that $G(X)$ contains at least one path from source to sink. The characteristic of $X_{FC}$ is: any vector $X$, if $X \ll X_{FC}$ (smaller than $X_{FC}$ according to BAT order), then $X$ must result in a disconnected network $G(X)$.

The identification of $X_{FC}$ is crucial for network reliability analysis because it can significantly reduce the number of vectors that need to be checked. By starting the search directly from $X_{FC}$ instead of from the zero vector, the algorithm can skip all vectors that are definitely disconnected, greatly improving computational efficiency.

The solution of $X_{FC}$ is based on the concept of minimum cut sets in graph theory in [24]. Formally, a subvector $X(1:k)$ constitutes a forward minimum cut if $G(X(1:k))$ is a disconnected graph and $G(X(1:k-1))$ is not disconnected, where $X(1:k)$ refers to the subvector in which arcs 1 through $k$ are set according to the values in $X$, and all other arcs are assumed to be functioning. This formulation precisely identifies the critical arc at position $k$ whose removal disconnects the network.

The $X_{FC}$ construction iteratively identifies the forward minimum cut set containing the arc with the smallest index, and sets this critical arc's state to 1. This process continues until the network becomes connected, ensuring that the resulting vector is indeed the lexicographically smallest connected vector.

**Algorithm: Find_$X_{FC}$**

**STEP 0:** Initialize all elements to 1 in $X_{FC}$, $k_0 = 0$, and set $i = 1$.

**STEP 1:** Identify the forward minimum cut set $X_{FC}(1: k_i)$ containing the arc with the smallest index $k_i$ in the current network configuration. If no such cut set exists (network is already connected), proceed to **STEP 3**.

**STEP 2:** Set all *j*-th coordinates in $X_{FC}$ to 0 for $k_{i-1} < j < k_i$, $i = i + 1$, then return to **STEP 1**.

**STEP 3:** Return $X_{FC}$ as the first connected vector.

The overall time complexity of the above algorithm is $O(|E|^3|V|)$, as it requires $O(|E|^2|V|)$ time for each iteration, with a maximum of *m* iterations. This is a significant improvement over the naive approach of checking all $2^m$ possible state vectors for connectivity.

To illustrate how the algorithm works, consider a simple network with 7 arcs as shown in Figure 3. Initially, $X_{FC} = (1, 1, 1, 1, 1, 1, 1)$ and the network is disconnected. In the first iteration, we identify the forward minimum cut set $(a_1, a_2) = (0, 0)$, set $X_{FC}(a_2) = 1$, and continue. After several iterations, the network becomes connected, and the algorithm terminates with $X_{FC} = (0, 1, 0, 0, 0, 1, 0)$ [24].

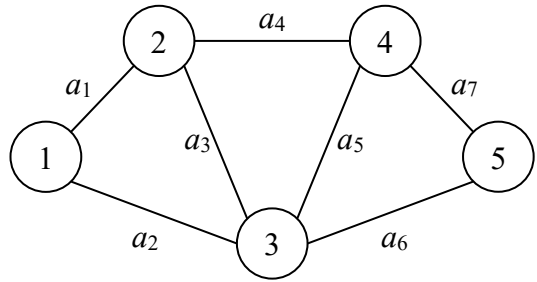

**Figure 3.** Example network for finding $X_{FC}$.

However, the above algorithm may fail in some cases. For example, in Figure 1, the correct $X_{FC}$ = (1, 0, 1, 0, 0, 0, 1, 0), but **Algorithm Find_$X_{FC}$** produces $X_{FC}$ = (0, 1, 0, 0, 1, 0, 0, 1). This discrepancy occurs because the algorithm identifies local minimum cuts in each iteration rather than considering the global optimal combination of arcs that would yield the lexicographically smallest connected vector. The algorithm's greedy approach works for many network topologies but can produce suboptimal results in networks with complex interdependencies among arcs [24].

It is worth noting that the concept of $X_{FC}$ is closely related to our proposed "earliest path" algorithm. While **Algorithm Find_$X_{FC}$** may not always identify the true first connected vector in the BAT ordering, the earliest path algorithm extends this concept by characterizing the specific path that first emerges in the BAT ordering. Our approach improves upon the $X_{FC}$ algorithm by using a more sophisticated binary weight allocation strategy that correctly identifies the earliest path even in complex network structures. Additionally, we provide faster computation methods for large-scale networks. The following sections will introduce our innovative approach that improves upon the $X_{FC}$ algorithm's efficiency while maintaining its theoretical guarantees.

## 5. EARLIEST AND LATEST PATH ALGORITHMS

In this section, we present the earliest and latest path algorithms that build upon the first connected vector concept and apply it to path optimization problems. By strategically assigning weights as powers of 2, we ensure that each combination of arcs corresponds to a unique total path weight, enabling efficient identification of paths with specific properties.

### 5.1 Problem Definition

The traditional shortest path problem is to find a path $P$ from source node 1 to sink node $n$ such that the sum of weights of all arcs on the path $\sum_{\forall a \in P} W(a)$ is minimized. However, in many practical application scenarios, in addition to path length, we also need to consider specific properties or priorities of path composition.

We define the earliest path problem as follows: identify the simple path $P$ from source node 1 to sink node $n$ that is discovered first in the BAT ordering sequence. This path represents the lexicographically smallest path that establishes connectivity between the source and sink nodes, corresponding to the first vector in the BAT sequence that contains a complete path.

Complementary to this, we define the latest path problem as follows: identify the simple path $P$ from source node 1 to sink node $n$ such that any vector found after this path in the BAT ordering sequence is guaranteed no simple paths.

By associating arc indices with specific attributes (such as reliability, priority, or risk), the earliest path problem naturally captures the needs of multi-objective path optimization while maintaining computational efficiency.

### 5.2 Binary Weight Allocation Mechanism

The key innovation in our approach is the binary weight allocation mechanism. For the earliest path algorithm, we assign a weight of $2^{(i-1)}$ to the $i$-th arc in the graph. For the latest path algorithm, we use the complementary weighting of $2^{(m-i)}$, where $m$ is the total number of arcs in the network. These weight allocations have the following characteristics:

- Each possible path combination has a unique weight sum
- The path's total weight directly maps to the binary representation of its included arc set
- For the earliest path, paths with smaller binary representations (using lower-indexed arcs) automatically receive lower total weights
- For the latest path, paths with larger binary representations (using higher-indexed arcs) automatically receive lower total weights

The main difference between these algorithms and the traditional shortest path is that: when there are multiple simple paths, our algorithms will choose the one with the smallest binary representation (earliest path) or the largest binary representation (latest path), which is particularly important for many practical applications, such as network reliability analysis, routing priority planning, and resource allocation optimization. Additionally, by adjusting the way arcs are numbered, the preference of path selection can be flexibly controlled, enabling the algorithm to adapt to specific needs in different domains.

### 5.3 Algorithm Design

Our path algorithms leverage Dijkstra's algorithm with the special binary weight allocations:

**Algorithm: Earliest path**

**Input:** Graph $G(V, E)$, source node 1, and sink node $n$

**Output:** The earliest path from 1 to $n$.

**STEP 0:** (Edge weight allocation) For each arc $a_i \in E$, set $W(a_i) = 2^{i-1}$

**STEP 1:** (Execute Dijkstra's algorithm)

    **STEP 1.1:** Initialize: For all nodes $v \in V$, set dist[$v$] = ∞, prev[$v$] = NULL, and dist[1] = 0.

    Create a priority queue Q containing all nodes, sorted by dist value

    **STEP 1.2:** While Q is not empty:

- $u$ = extract node with smallest dist value from Q

- if $u = n$: break the loop

- for each arc $e_{u,v} \in E$: alt = dist[$u$] + $W(e_{u,v})$, if alt < dist[$v$], set dist[$v$] = alt, prev[$v$] = $u$, and update $v$'s priority in Q

**STEP 2:** (Reconstruct the earliest path):

    **STEP 2.1:**    Set path = empty list, current = $n$

    **STEP 2.2:**    while current ≠ NULL:

- insert current at the front of path

- set current = prev[current];

    **STEP 2.3:**    return path

**Algorithm: Latest path**

**Input:**    Graph $G(V, E)$, source node 1, and sink node $n$, and total number of arcs $m$

**Output:**    The latest path from 1 to $n$.

**STEP 0:** (Edge weight allocation) For each arc $a_i \in E$, set $W(a_i) = 2^{(m-i)}$

**STEP 1:** (Execute Dijkstra's algorithm) - Same as the earliest path algorithm

**STEP 2:** (Reconstruct the path) - Same as the earliest path algorithm

Due to the special properties of our weight allocation, the paths found by these algorithms are not only valid paths from source to sink but also the ones with the smallest (for earliest) or largest (for latest) binary representation among all valid paths. This elegantly combines the BAT ordering concept with Dijkstra's algorithm for efficient path discovery.

### 5.4 Time Complexity Analysis

The time complexity of our algorithms can be analyzed step by step:

**STEP 0:** Arc weight allocation requires $O(|E|)$ time, traversing each arc and setting its weight.

**STEP 1:** Dijkstra's algorithm implemented with a binary heap has a time complexity of $O((|V|+|E|)\log|V|)$. This complexity derives from several operations: initialization requires $O(|V|)$ time; building the priority queue requires $O(|V|\log|V|)$ time; the main loop executes

|V| times with each extraction operation requiring $O(\log|V|)$ time; and each arc is processed at most once, with each priority queue update requiring $O(\log|V|)$ time.

**STEP 2:** Path reconstruction requires $O(|V|)$ time, as the path contains at most $|V|$ nodes.

Therefore, the overall time complexity of both algorithm is: $O(|E|) + O((|V|+|E|)\log|V|) + O(|V|) = O((|V|+|E|)\log|V|)$. This is the same as the time complexity of the standard Dijkstra algorithm, indicating that our path algorithms are computationally equivalent to the traditional Dijkstra algorithm, with only a linear-time preprocessing step (edge weight allocation) added.

## 5.5 Proof of Correctness for the Earliest Path Algorithm

To prove the correctness of our algorithms, we need to establish two key properties for each:

1. The algorithm finds a valid path from source to sink if one exists
2. The path found is optimal according to our specific optimization criteria (earliest or latest in the BAT ordering)

### 5.5.1 Earliest Path Algorithm

The earliest path algorithm employs a binary weight allocation strategy where lower-indexed arcs receive exponentially smaller weights. This strategic weighting creates a direct correspondence between the BAT ordering of paths and their total weight, allowing Dijkstra's algorithm to naturally identify the lexicographically smallest path.

**Theorem 1:** The proposed earilest algorithm finds a valid path from source to sink.

**Proof:** This follows directly from the correctness of Dijkstra's algorithm. Since we are using Dijkstra's algorithm with non-negative weights ($2^{(i-1)} \geq 0$ for all $i$), the algorithm correctly finds a path from source to sink if one exists, minimizing the sum of weights along the path.

**Theorem 2:** The path found is the earliest in the BAT ordering.

**Proof:** Let's proceed by contradiction. Suppose $P = (p_1, p_2, \ldots, p_m)$ is the path found by our algorithm, and there exists another path $X = (x_1, x_2, \ldots, x_m)$ that is earlier in the BAT ordering, i.e., $x_i < p_i$ and $x_j = p_j$ for all $j > i$. By definition of the BAT ordering, if $X$ is earlier than P, then there exists an index i such that:

- P includes arc $a_i$
- X does not include arc $a_i$
- For all arcs $a_j$ with $i < j$, either both paths exclude $a_j$

Given our weight allocation where $W(a_k) = 2^{(k-1)}$, we can express the weights as

$$W(P) = \sum_{\forall P(a)=1} W(a) = \sum_{k=1}^{i} 2^{k-1} p_k \tag{1}$$

$$W(X) = \sum_{\forall X(a)=1} W(a) = \sum_{k=1}^{i} 2^{k-1} x_k. \tag{2}$$

Beause $x_i < p_i$ and

$$2^{i-1} > \sum_{k=1}^{i-1} 2^{k-1} = \frac{2^{i-1}-1}{2-1} = 2^{i-1} - 1, \tag{3}$$

we have $W(P) > W(X)$, which contradicts our assumption that P is the shortest path found by Dijkstra's algorithm (which minimizes total weight). Hence, there cannot exist a path X that is earlier in the BAT ordering than the path P found by our algorithm. This proves that our algorithm correctly identifies the earliest path.

### 5.5.2 Latest Path Algorithm

For the latest path algorithm, we use the complementary weight assignment where higher-indexed arcs receive smaller weights, effectively reversing the priorities compared to the earliest path algorithm. This approach ensures that paths containing predominantly higher-indexed arcs are preferred, leading to the identification of the lexicographically latest path.

**Theorem 3:** The proposed latest algorithm finds a valid path from source to sink.

**Proof:** This follows directly from the correctness of Dijkstra's algorithm with non-negative weights ($2^{(m-i)} > 0$ for all $i$).

**Theorem 4:** The path found is the latest in the BAT ordering.

**Proof:** Using similar logic to the earliest path proof, but with the complementary weight allocation $W(a_i) = 2^{(m-i)}$, we can show that the path with the largest binary representation (latest in BAT ordering) will have the smallest weight sum under this allocation. This ensures that Dijkstra's algorithm finds the latest path according to the BAT ordering.

**Corollary 1.** (Optimality with respect to traditional distance) If we consider two paths with the same number of arcs or the same traditional distance, our weight allocation ensures that the path with the smallest binary representation (i.e., the earliest in BAT ordering) or the largest binary representation (i.e., the latest in BAT ordering) will have the smallest weight sum in their respective algorithms.

**Proof:** This follows from the property that for any set of arcs $\{a_i, a_j, \ldots, a_k\}$ and another set $(a_p, a_q, \ldots, q_r)$ where the first set precedes the second in BAT ordering, the sum of weights

$$2^{(i-1)} + 2^{(j-1)} + \ldots + 2^{(k-1)} < 2^{(p-1)} + 2^{(q-1)} + \ldots + 2^{(r-1)} \tag{3}$$

for the earliest path algorithm. For the latest path algorithm with weights $2^{(m-i)}$, the inequality is reversed, ensuring the latest path is found.

Therefore, our algorithms not only find valid paths but specifically find the earliest and latest paths according to the BAT ordering, completing our proof of correctness.

## 6. EXAMPLES AND PRACTICAL APPLICATIONS

This section demonstrates the application of the earliest and latest path algorithms through specific examples and discusses their value in practical scenarios. We first use a simple network topology to illustrate the execution process and results of the algorithms, clearly showing how binary weight allocation affects path selection. Subsequently, we discuss real-world applications of the algorithms in complex network environments, demonstrating their advantages in solving multi-objective path optimization problems.

### 6.1 Examples

As network scale increases, computational costs grow exponentially. Figure 1 shows a widely used benchmark network in binary state networks, with a scale suitable for detailed explanation of our proposed algorithm. Suppose we have a network as shown in Figure 1:

Nodes: 1, 2, 3, 4, 5, 6 (where 1 is the source node and 6 is the target node)

Arcs: $a_1 = e_{1,2}$, $a_2 = e_{1,3}$, $a_3 = e_{2,4}$, $a_4 = e_{2,5}$, $a_5 = e_{3,5}$, $a_6 = e_{4,5}$, $a_7 = e_{4,6}$, $a_8 = e_{5,6}$

### 6.1.1 Earliest Path Example

Now we execute the algorithm:

**STEP 0:** According to the earliest path algorithm, we first assign weights to each edge as shown in Table 2:

Table 2. Arcs and weights for earliest path

| $i$ | $W(a_i)$ | $i$ | $W(a_i)$ |
|---|---|---|---|
| 1 | $2^0 = 1$ | 5 | $2^4 = 16$ |
| 2 | $2^1 = 2$ | 6 | $2^5 = 32$ |
| 3 | $2^2 = 4$ | 7 | $2^6 = 64$ |
| 4 | $2^3 = 8$ | 8 | $2^7 = 128$ |

**STEP 1:** Execute Dijkstra's main loop

    **STEP 1.1:** (Initialization) Set the distance of all nodes to infinity: dist[1, 2, ..., 6] = [0, ∞, ∞, ∞, ∞, ∞], the predecessor of all nodes to NULL: prev[1, 2, ..., 6] = [NULL, NULL, NULL, NULL, NULL, NULL], and set unvisited node Q = {1, 2, 3, 4, 5, 6}.

    **STEP 1.1:** Because Q is not empty:

        **Iteration 1:** Select node 0 (distance is 0). Update node 1: dist[1] = min(∞, 0+1) = 1, prev[1] = 0. Update node 2: dist[2] = min(∞, 0+2) = 2, prev[2] = 0. Remove node 0: Q = {1, 2, 3, 4, 5}.

        **Iteration 2:** Select node 1 (distance is 1). Update node 3: dist[3] = min(∞, 1+4) = 5, and prev[3] = 1. Update node 4: dist[4] = min(∞, 1+8) = 9, and prev[4] = 1. Remove node 1: Q = {2, 3, 4, 5}.

        **Iteration 3:** Select node 2 (distance is 2). Update node 4: dist[4] = min(9, 2+16) = 9, and prev[4] remains unchanged. Remove node 2: Q = {3, 4, 5}.

        **Iteration 4:** Select node 3 (distance is 5). Update node 5: dist[5] = min(∞, 5+64) = 69, prev[5] = 3. Remove node 3: Q = {4, 5}.

        **Iteration 5:** Select node 4 (distance is 9). Update node 5: dist[5] = min(69, 9+128) = 69, and prev[5] remains unchanged. Remove node 4: Q = {5}.

**Iteration 6:** Select node 5 (distance is 37). Because Q is empty, algorithm ends.

**STEP 2:** (Reconstruct the earliest path): Start backtracking from node 5: path = [5, 3, 1, 0]. After reversal, get the final path: 0 → 1 → 3 → 5.

This path 0 → 1 → 3 → 5 is the earliest path among all paths from node 0 to node 5. Alternative paths exist, such as 0→2→4→5 with a weight sum of 2+16+128=146 and 0→1→4→5 with a weight sum of 1+8+128=137. Although all three are valid paths with the same number of arcs (three), our algorithm selects 0→1→3→5 because it has the smallest binary weight sum, corresponding to the smallest binary representation in the BAT ordering.

### 6.1.2 Latest Path Example

Now we execute the latest path algorithm on the same network, but with complementary weights as shown in Table 3:

Table 3. Arcs and weights for latest path.

| $i$ | $W(a_i)$ | $i$ | $W(a_i)$ |
|---|---|---|---|
| 1 | $2^7 = 128$ | 5 | $2^3 = 8$ |
| 2 | $2^6 = 64$ | 6 | $2^2 = 4$ |
| 3 | $2^5 = 32$ | 7 | $2^1 = 2$ |
| 4 | $2^4 = 16$ | 8 | $2^0 = 1$ |

**STEP 1:** Execute Dijkstra's main loop (similar initialization and iteration process as before)

After running the algorithm with these weights, we obtain the latest path: 1 → 3 → 5 → 6.

This path has a weight sum of 64 + 8 + 1 = 73 under the latest path weighting scheme. Despite having the same number of arcs as the earliest path, this path has the highest indices possible, making it the lexicographically latest path in the BAT ordering. Any vector that comes after this path in the BAT ordering is guaranteed to contain no simple path from source to sink.

The combination of the earliest path (1→2→4→6) and latest path (1→3→5→6) effectively partitions the solution space into three regions:

- Vectors before 1→2→4→6: guaranteed disconnection

- Vectors between 1→2→4→6 and 1→3→5→6: transitional connectivity

- Vectors after 1→3→5→6: guaranteed no simple paths

This partitioning provides significant computational benefits for network reliability analysis and similar applications.

### 6.2 Practical Applications

The earliest and latest path algorithms based on binary weights demonstrate significant application value in multiple domains, effectively solving multi-objective optimization problems in complex networks through their unique weight allocation mechanisms:

#### 6.2.1 Network Reliability Analysis

The algorithms naturally guide path selection towards more reliable component combinations by mapping less critical network elements to smaller index values. In telecommunications network reliability assessment, these methods can identify vulnerable links that are difficult to discover using traditional methods, providing more accurate failure probability estimates.

More importantly, by identifying both the earliest and latest paths, these algorithms establish crucial boundaries for network connectivity. Any configuration with arcs below the earliest path threshold will definitely be disconnected, while any configuration with arcs above the latest path threshold will definitely be connected. This allows network designers to focus computational resources on the transitional region between these thresholds, greatly enhancing analysis efficiency for large networks.

By iteratively applying the algorithms and removing arcs from discovered paths, multiple disjoint paths can be systematically identified, establishing a more comprehensive network redundancy structure. This approach enhances the overall reliability and interference resistance of the system by ensuring diverse routing options that avoid common points of failure.

#### 6.2.2 Communication System Design

In 5G and future 6G network planning, the algorithms provide an ingenious method to balance transmission efficiency and network security. By assigning index values reflecting priorities to different types of communication links (such as fiber optic, microwave, satellite links), the algorithms optimize path composition while meeting bandwidth and latency requirements.

When designing backup routing schemes, the complementary nature of the earliest and latest paths provides naturally diverse routing options. The earliest path typically utilizes the most reliable or preferred links, while the latest path offers a maximally diverse alternative route, providing the network with more efficient fault recovery mechanisms and reducing service interruption time during network failures.

### 6.2.3 Traffic Planning and Smart Transportation Systems

In urban traffic network optimization, the algorithms generate optimal routes that balance transit efficiency with other decision factors by assigning larger index values to congested road segments or environmentally sensitive areas. This approach reduces average commute times while simultaneously decreasing traffic volume through residential areas.

The algorithms also excel in dynamic routing scenarios, where real-time traffic conditions can be incorporated by adjusting arc indices based on current congestion levels. This adaptive approach outperforms traditional routing algorithms during peak traffic hours, providing more resilient transportation networks that can effectively respond to unexpected disruptions like accidents or construction.

### 6.2.4 Supply Chain Optimization

In global supply chain management, the earliest and latest path algorithms offer significant advantages for optimizing multi-modal transportation routes. By encoding factors such as transportation cost, transit time, customs inspection frequency, and environmental impact into arc indices, the algorithms can identify routes that balance these competing objectives according to business priorities.

The binary weight allocation proves particularly valuable for handling complex trade-offs between cost-efficiency and customs processing time, a critical factor in cross-border logistics that traditional shortest path algorithms struggle to adequately address.

The above practical applications demonstrate that the earliest and latest path algorithms not only offer theoretical elegance through their binary weight allocation mechanisms but also deliver substantial performance improvements in real-world scenarios. Their ability to implicitly encode multiple objectives without additional computational overhead makes them invaluable tools for network optimization across diverse domains.

## 7. CONCLUSION AND FUTURE PROSPECTS

This paper presents earliest and latest path algorithms based on the Dijkstra algorithm, which effectively find paths with optimal binary representations by assigning exponential weights of $2^{(i-1)}$ and $2^{(m-i)}$ to the $i$-th arc in the network, respectively. Our approach leverages the Binary Addition Tree (BAT) ordering and first connected vector concepts to identify lexicographically smallest and largest paths that establish connectivity between source and sink nodes.

The key theoretical contributions of this work include: (1) a formal definition of earliest and latest paths within the BAT framework, (2) a binary weight allocation mechanism that implicitly encodes path priorities, and (3) a computationally efficient algorithm that maintains the same time complexity as the traditional Dijkstra algorithm while providing multi-objective optimization capabilities. The examples and case studies presented in Section 6 demonstrate the algorithm's effectiveness across diverse domains including telecommunications, transportation networks, and supply chain management.

Compared to traditional shortest path algorithms, our approach partitions the solution space into three regions (guaranteed disconnection, transitional connectivity, and guaranteed no simple paths), allowing for more efficient exploration of the solution space. The earliest path represents the first vector in the BAT ordering that contains a complete path, while the latest path ensures that any vector found after it will have no simple paths. By identifying both earliest and latest paths, our method provides important computational boundaries that can significantly reduce the number of vectors that need to be evaluated in large-scale networks.

While our approach offers significant advantages, it does have limitations. The current implementation assumes static networks with fixed arc indices, and the binary weight allocation may lead to numerical precision issues for very large networks due to the exponential growth of weights.

Several promising directions for future research include:

1. **Dynamic Network Adaptability and Self-Learning Mechanisms**: Investigating real-time optimization strategies for networks with frequently changing topologies or arc weights, potentially integrating machine learning techniques to automatically determine optimal arc index configurations based on historical performance data.

2. **Extended Problem Formulations**: Adapting the algorithm to handle multi-source multi-sink problems, multi-commodity flow scenarios, and constrained shortest path problems to address more complex network applications.

3. **Integration with Advanced Optimization Techniques**: Exploring combinations with metaheuristic approaches such as genetic algorithms or ant colony optimization to enhance performance for ultra-large-scale networks where exact methods become computationally prohibitive.

4. **Theoretical Extensions**: Further developing the mathematical foundations of the relationship between BAT ordering and network reliability, potentially establishing formal bounds on reliability metrics based on earliest and latest path properties.

5. **Expanded Application Domains**: Investigating the algorithm's application in emerging fields such as quantum network routing, blockchain optimization, and complex biological network analysis, where multi-objective path optimization plays a crucial role.

These research directions will not only enrich the theoretical foundation of our approach but also significantly expand its practical application value. By bridging the gap between elegant mathematical formulations and real-world network optimization challenges, our work provides a framework for developing more comprehensive and effective solutions for complex systems across multiple domains.

## ACKNOWLEDGEMENTS

This research was supported in part by the Ministry of Science and Technology, R.O.C. under grant MOST 107-2221-E-007-072-MY3 and MOST 110-2221-E-007-107-MY3. This article was once submitted to arXiv as a temporary submission that was just for reference and did not provide the copyright.

**REFERENCES**

[1] WC Yeh, JS Lin (2018) New parallel swarm algorithm for smart sensor systems redundancy allocation problems in the Internet of Things. The Journal of Supercomputing 74, 4358-4384.

[2] D Narciandi-Rodriguez, J Aveleira-Mata, M. T. García-Ordás, J. Alfonso-Cendón, C. Benavides, H. Alaiz-Moretón (2025) A cybersecurity review in IoT 5G networks. Internet of Things 30, 101478.

[3] WC Yeh (2009) A simple universal generating function method to search for all minimal paths in networks. IEEE Transactions on Systems, Man, and Cybernetics-Part A: Systems and Humans 39 (6), 1247-1254.

[4] WC Yeh, W Zhu, CL Huang, TY Hsu, Z Liu, SY Tan (2022) A New BAT and PageRank Algorithm for Propagation Probability in Social Networks. Applied Sciences 12, doi.org/10.3390/app12146858.

[5] WC Yeh, W Zhu, CL Huang (2022) Reliability of Social Networks on Activity-on-Node Binary-State with Uncertainty Environments. Applied Sciences 2022 (12), doi:10.3390/app12199514.

[6] Z Liu, WC Yeh, KY Lin, CSH Lin, CY Chang (2024) Machine learning based approach for exploring online shopping behavior and preferences with eye tracking. Computer Science and Information Systems 21 (2), 593-623.

[7] WC Yeh, SC Wei (2012) Economic-based resource allocation for reliable Grid-computing service based on Grid Bank. Future Generation Computer Systems 28 (7), 989-1002.

[8] WC Yeh, W Zhu, Y Yin, CL Huang (2023) Cloud Computing Considering Both Energy and Time Solved by Two-Objective Simplified Swarm Optimization. Applied Sciences 13 (2077), 10.3390/app13042077.

[9] J Wang, WC Yeh, NN Xiong, J Wang, X He, CL Huang (2019) Building an improved Internet of Things smart sensor network based on a three-phase methodology. IEEE Access 7, 141728-141737.

[10] C Luo, B Sun, K Yang, T Lu, WC Yeh (2019) Thermal infrared and visible sequences fusion tracking based on a hybrid tracking framework with adaptive weighting scheme. Infrared Physics & Technology 99, 265-276.

[11] PC Su, SY Tan, ZY Liu, WC Yeh (2022) A Mixed-Heuristic Quantum-Inspired Simplified Swarm Optimization Algorithm for scheduling of real-time tasks in the multiprocessor system. Applied Soft Computing, 10.1016/j.asoc.2022.109807.

[12] EW Dijkstra (1959) A Note on Two Problems in Connexion with Graphs. Numerische Mathematik 1(1), 269-271.

[13] R Bellman (1958) On a Routing Problem. Quarterly of Applied Mathematics 16(1), 87-90.

[14] RW Floyd (1962) Algorithm 97: Shortest Path. Communications of the ACM 5(6), 345.

[15] WC Yeh, SY Tan, M Forghani-elahabad, M ElKhadiri, Y Jiang, CS Lin (2022) New Binary-Addition Tree Algorithm for the All-Multiterminal Binary-State Network Reliability Problem. Reliability Engineering and System Safety 224, doi.org/10.1016/j.ress.2022.108557.

[16] Z Hao, WC Yeh (2025) GE-MBAT: An Efficient Algorithm for Reliability Assessment in Multi-State Flow Networks. Reliability Engineering & System Safety, 110916.

[17] WC Yeh (2025) Enhancing Reliability Calculation for One-Output k-out-of-n Binary-state Networks Using a New BAT. Reliability Engineering and System Safety, 10.1016/j.ress.2025.110835.

[18] WC Yeh, M Forghani-elahabad (2024) An efficient parallel approach for binary-state network reliability problems. Annals of Operations Research, 1-22.

[19] WC Yeh, W Zhu (2024) Optimal Allocation of Financial Resources for Ensuring Reliable Resilience in Binary-State Network Infrastructure. Reliability Engineering & System Safety 250, 110265.

[20] WC Yeh (2024) Time-reliability optimization for the stochastic traveling salesman problem. Reliability Engineering & System Safety 248, 110179.

[21] WC Yeh (2024) A New Hybrid Inequality BAT for Comprehensive All-Level d-MP Identification Using Minimal Paths in Multistate Flow Network Reliability Analysis. Reliability Engineering & System Safety 244, 109876.

[22] WC Yeh, CC Kuo (2020) Predicting and modeling wildfire propagation areas with BAT and maximum-state PageRank. Applied Sciences 10 (23), 8349.

[23] WC Yeh, E Lin, CL Huang (2021) Predicting Spread Probability of Learning-Effect Computer Virus. Complexity 2021, 6672630.

[24] WC Yeh (2021) A Quick BAT for Evaluating the Reliability of Binary-State Networks. Reliability Engineering & System Safety 216, 107917.

[25] Z Hao, WC Yeh (2025) Applying Incremental Learning in Binary-Addition-Tree Algorithm in Reliability Analysis of Dynamic Binary-State Networks. Reliability Engineering & System Safety 261, 111072.

[26] ML Fredman, RE Tarjan (1987) Fibonacci heaps and their uses in improved network optimization algorithms. Journal of the ACM (JACM), 34(3), 596–615.

[27] WC Yeh (2022) BAT-based Algorithm for Finding All Pareto Solutions of the Series-Parallel Redundancy Allocation Problem with Mixed Components. Reliability Engineering & System Safety 228, 10.1016/j.ress.2022.108795.

[28] PE Hart, NJ Nilsson, B Raphael (1968) A formal basis for the heuristic determination of minimum cost paths. IEEE Transactions on Systems Science and Cybernetics, 4(2), 100–107.